\documentclass[11pt,twoside]{article}
\usepackage{mathbbold}
\usepackage{amsmath, amssymb, mathrsfs, graphicx}

\def\scr{\mathscr}

\def\az{\alpha}  \def\bz{\beta}
    \def\dz{\delta}
    \def\fz{\varphi}
\def\gz{\gamma}  \def\kz{\kappa}
\def\lz{\lambda}

\def\qz{\psi}
        
        \def\uz{\theta}
\def\vz{\varepsilon}

\def\qd{\quad}
\def\qqd{\qquad}

\newcommand{\mathsym}[1]{{}}

\def\scr{\mathscr}

\def\le{\leqslant}
\def\ge{\geqslant}

\font\cms=cmss9 scaled \magstep1

\def\nnd{\noindent}

\def\thm{\nnd\bg{thm1}}

\def\prp{\nnd\bg{prp1}}

\def\xmp{\nnd\bg{xmp1}}

\def\dethm{\end{thm1}}

\def\deprp{\end{prp1}}
\def\dexmp{\end{xmp1}}

\def\bg{\begin}
\def\be{\bg{equation}}
\def\de{\end{equation}}

\def\dear{\end{eqnarray}}
\def\lb{\label}
\def\ct{\cite}

\newcommand{\rf}[2]{[\ref{#1}; #2]}

\def\den{\end{enumerate}}

\def\d{\text{\rm d}}

\begin{document}

\allowdisplaybreaks[4]
\thispagestyle{empty}
\renewcommand{\thefootnote}{\fnsymbol{footnote}}

\noindent {\small Festschrift Masatoshi Fukushima:
In Honor of Masatoshi Fukushima's Sanju}\newline
{\small Eds: Z.Q. Chen, N. Jacob, M. Takeda, T. Uemura.}\newline
{\small World Scientific, 2015}

\vspace*{.5in}
\begin{center}
{\bf\Large Progress on Hardy-type Inequalities}
\vskip.15in {Mu-Fa Chen}
\end{center}
\begin{center} (Beijing Normal University, Beijing 100875, China)\\
\vskip.1in
October 4, 2013
\end{center}
\vskip.1in

\markboth{\sc Mu-Fa Chen}{\sc Progress on Hardy-type inequalities}

\footnotetext{2000 {\it Mathematics Subject Classifications}.\quad 26D10, 60J60, 34L15.}
\footnotetext{{\it Key words and phrases}.\quad
Hardy-type inequality, stability speed, optimal constant, basic estimate.}

\bigskip

\begin{abstract} This paper surveys some of our recent progress on
Hardy-type inequa\-lities which consist of a well-known topic in Harmonic Analysis. In the first section, we recall the original probabilistic motivation dealing with the stability speed in terms of the $L^2$-theory.
A crucial application of a result by Fukushima and Uemura (2003)
is included.
In the second section, the non-linear case (a general Hardy-type inequality) is handled with a direct and analytic proof. In the last section, it is illustrated that the basic estimates presented in the first two sections can still be improved considerably.
\end{abstract}

\medskip

This paper mainly concerns with the following Hardy-type inequality
\be\bigg(\int_{-M}^N |f-\pi(f)|^q\d\mu\bigg)^{1/q}\le A
 \bigg(\int_{-M}^N \big|f'\big|^p\d\nu\bigg)^{1/p},\lb{01}\de
where $p, q \in (1, \infty)$, $\mu$ and $\nu$ are Borel measures
on an interval $[M, N]\,(M, N\le\infty)$.
Here, we assume that $\mu[-M, N]<\infty$ and
define a probability measure $\pi=(\mu [-M, N])^{-1}\mu$.
Then $\pi(f):=\int f \d\pi$. The functions $f$ are assumed to be
absolutely continuous on $(-M, N)$ and belong to $L^q(\mu)$.
For simplicity, we may also write the inequality as
$$\|f-\pi(f)\|_{L^q(\mu)}\le A\, \big\|f'\big\|_{L^p(\nu)}.$$
To save our notation, assume the constant $A$ to be optimal.
The linear case that $p=q=2$ is discussed in the next section, where
a result by Fukushima and Uemura \ct{fmut} plays an important role. The general
case is studied in Section 2. In the last section, we show the possibility
for improving further the basic estimates of the optimal constant.

\section{Linear case: $p=q=2$.}

Let us recall the original probabilistic problem.
Throughout this section, we fix $p=q=2$.
Consider a second-order elliptic operator on $(-M, N)$:
$$L= a(x)\frac{\d^2}{\d x^2}+ b(x)\frac{\d}{\d x}, \qqd a(x)>0 \text{ on } (-M, N).$$
Then the two measures used in inequality (\ref{01}) are as follows
\be\mu(\d x)=\frac{e^{C(x)}}{a(x)}\d x, \qqd \nu(\d x)= e^{C(x)}\d x,\qqd
C(x):=\int_{\uz}^x \frac{b}{a},\de
here in the last integral and in what follows, the Lebesgue measure $\d x$ is omitted, $\uz\in (-M, N)$ is a reference point.
Denote by $\{P_t\}_{t\ge 0}$ the (maximal) semigroup generated by $L$
on $L^2(\mu)$. Here ``maximal'' means the Dirichlet form having the maximal
domain which we learnt earlier from Fukushima \ct{fm80} (cf. \rf{cmf04}{\S 6.7}).
We are interested in the stability speed, for instance, the $L^2$-exponential
convergence rate $\vz$:
$$\|P_t f-\pi(f)\|_{L^2(\mu)}\le \|f-\pi(f)\|_{L^2(\mu)}e^{-\vz t}, \qqd t\ge 0.$$
Then, it turns out that the largest rate $\vz_{\max}$ coincides with $A^{-1}$ given in (\ref{01}) (cf. \rf{cmf04}{Theorem 9.1}).

We can state one of our recent results as follows.

\thm\lb{t1-1}{\rm\rf{cmf10}{Theorem\;10.2}}\; {\cms Let $a>0$, $a$ and $b$ be continuous on $[-M,
N]$ (or $(-M, N]$ if $M=\infty$, for instance). Assume that $\mu(-M,
N)<\infty$. Then for the optimal constant $A$, we have the basic estimates:
${\kz}\le A\le 2 {\kz}$, where
$$\kz^{-2}=\inf_{-M<x<y<N}\big[\mu(-M, x)^{-1}+\mu(y, N)^{-1}\big]{\hat\nu}(x, y)^{-1},$$
$\mu (x, y)=\int_x^y \d\mu$, and ${\hat\nu}(\d x)=e^{- C(x)}\d x$.}
\dethm

The continuity assumption on $a$ and $b$ is not essential and will be removed in the next
section. To understand the proof of this theorem, assume that $M, N<\infty$. Then the general
case can be done by an approximating procedure. The optimal constant actually describes an
eigenvalue $\lz_1(=A^{-2})$ defined by
$$\lz_1=\inf\Big\{\big\|f'\big\|_{L^2(\nu)}^2\!: \pi(f)=0,\; \|f\|_{L^2(\mu)}=1\Big\}.$$
Let $g$ be the eigenfunction corresponding to $\lz_1$:
$$L g = - \lz_1 g, \qqd g\ne 0.$$
To ignore a constant (say $\pi(g)$) from $g$, making derivative on both sides of
the equation and replacing $g'$ by $f$, we get
$$L_S f=-\lz_1 f,$$
where
$$L_S=a(x)\frac{\d^2}{\d x^2} + \big(a'(x)+b(x)\big)\frac{\d}{\d x}+ b'(x)$$
which is a Schr\"odinger operator. Because $g'(-M)=0=g'(N)$, the boundary condition for
$L_S$ becomes $f(-M)=0=f(N)$. This leads to the principal eigenvalue of $L_S$:
$$\lz_S=\sup_{f\in {\scr F}}\inf_{x\in (-M, N)}\frac{-L_S f}{f}(x),$$
where
$${\scr F}=\big\{f\in {\scr C}^1[-M, N]\cap {\scr C}^2(-M, N): f(-M)=f(N)=0,\; f|_{(-M, N)}>0\big\}.$$
Note that the zero-point of $g$ in the original eigenequation for $L$ is located inside of
the interval $(-M, N)$, not explicitly known; the zero-points for the eigenvalue $\lz_S$ are located
at the boundaries $-M$ and $N$ only. This is the advantage of $L_S$. However, there is an extra
term $b'$ in operator $L_S$ which costs some trouble as usual. To avoid this, we rewrite $L$ as
$$L=\frac{\d}{\d\mu}\frac{\d}{\d{\hat\nu}}. $$
Then we define a dual operator of $L$ as follows.
$$L^*=\frac{\d}{\d\mu^*}\frac{\d}{\d{\hat\nu}^*}
:=\frac{\d}{\d{\hat\nu}}\frac{\d}{\d\mu},$$
i.e. an exchange of $\mu$ and $\hat\nu$. More explicitly,
$$L^*=a(x)\frac{\d^2}{\d x^2} + \big(a'(x)-b(x)\big)\frac{\d}{\d x}.$$
Next, define
$$\lz_0^*=\sup_{f^*\in {\scr F}}\inf_{x\in (-M, N)}\frac{-L^* f^*}{f^*}(x).$$
In view of the next result which is crucial in proving Theorem \ref{t1-1},
it is clear that we have thus removed the extra term $b'$ in the operator $L_S$.

\prp \qd {\cms We have $A^{-2}=\lz_1=\lz_S=\lz_0^*$.}
\deprp

Here, in proving $\lz_1=\lz_S$, we have used a mathematical tool --- the coupling technique
(cf. \ct{cmwf97} or \ct{cmf05}). We have also used another tool --- dual technique in proving $\lz_S=\lz_0^*$. It is interesting that
they are the main tools used in the study on interacting particle systems (cf. \ct{cmf04} and
references therein). To obtain the basic estimates listed in Theorem \ref{t1-1}, we need one more mathematical tool --- the capacitary method.
The next result is taken from Fukushima \& Uemura \ct{fmut} and \ct{cmf05, cmf05b}, see also \ct{fot11}.

\thm\lb{t1-3} \qd{\cms For a regular transient Dirichlet form  $(D, {\scr D}(D))$
with locally compact state space $(E, {\scr E})$,
the optimal constant $A_{\mathbb B}$ in the Poincar{\cms\'e}-type
inequality
$$\big\|f^2\big\|_{\mathbb B}\le A_{\mathbb B}^2\, D(f),\qqd f\in {\scr C}_K^\infty (E)$$
satisfies
$ B_{\mathbb B}\le A_{\mathbb B}\le 2 B_{\mathbb B}, $
where $\|\cdot\|_{\mathbb B}$ is the norm in a normed linear space ${\mathbb B}$ and}
$$B_{\mathbb B}^2=\sup_{\text{\rm compact}\, K} {\text{\rm Cap}(K)}^{-1}{{\|\mathbbold{1}_K\|_{\mathbb B}}}. $$
\dethm

The space ${\mathbb B}$ can be very general, for instance $L^p(\mu)\,(p\ge 1)$ or the
Orlicz spaces. In the present context,
$D(f)= \int_{-M}^N {f'}^2 e^C=\|f'\|_{L^2(\nu)}^2$,
${\scr D}(D)$ is the closure of ${\scr C}_K^{\infty}(-M, N)$
with respect to the norm $\|\cdot\|_D$: $\|f\|_D^2\!=\!\|f\|^2\!+\!D(f)$, and
$$\text{\rm Cap}(K)=\inf\big\{D(f): f\in {\scr C}_K^{\infty}(-M, N),\, f|_K\ge 1\big\}.$$

Note that we have the
universal factor 2 here and the isoperimetric constant
$B_{\mathbb B}$ has a very compact form.
We now need to compute the capacity only. The problem is that the capacity
is usually not  explicitly computable. For instance, at the moment, we do not
know how to compute it for Schr\"odinger operators even for the
elliptic operators having killings. Very lucky,
we are able to compute the capacity for the one-dimensional
elliptic operators. The result has a simple expression:
$$B_{{\mathbb B}}^2=\sup_{-M< x < y <N}
\big[{\hat\nu}(-M, x)^{-1} + {\hat\nu}(y, N)^{-1}\big]^{-1}\|\mathbbold{1}_{(x,\, y)}\|_{{\mathbb B}}.$$
It looks strange to have double inverse here. So, making inverse in both sides, we get
$$B_{{\mathbb B}}^{-2}=\inf_{-M< x < y <N}
\big[{\hat\nu}(-M, x)^{-1} + {\hat\nu}(y, N)^{-1}\big]\,\|\mathbbold{1}_{(x, y)}\|_{{\mathbb
B}}^{-1}.$$
Applying this result to ${\mathbb B}=L^1(\mu)$, we obtain
the solution to the case having double Dirichlet boundaries: $\lz_0=A_{L^1(\mu)}^{-2}$ and
$$\kz_0^{-2}=B_{L^1(\mu)}^{-2}=\inf_{-M< x < y <N}\big[{\hat\nu}(-M, x)^{-1} + {\hat\nu}(y, N)^{-1}\big]\,
{\mu (x, y)^{-1}}.$$

Applying the last result to the dual process, we have
not only
$$\big(\kz_0^*\big)^2\le \lz_1=\lz_S=\lz_0^*\le 4 \big(\kz_0^*\big)^2$$
but also
$$\aligned
{\big(\kz_0^{*}\big)^{-2}}&=
\inf_{x < y}\big[{\hat\nu}^*(-M, x)^{-1} + {\hat\nu}^*(y, N)^{-1}\big]\,
\mu^* (x, y)^{-1}\\
&=\inf_{x < y}\big[\mu(-M, x)^{-1} + \mu(y, N)^{-1}\big]\,
{\hat\nu}(x, y)^{-1}\\
&=\kz^{-2}.
\endaligned$$
We have thus arrived at the assertion of Theorem \ref{t1-1}. Refer to \rf{cmf10}{\S 10} and \ct{cmf12} for more details.

To conclude this section, we remark that the use of the capacity
is natural in the higher dimensions, since in which the boun\-dary may be very complicated. However, it seems unnecessary to use it in the present
one-dimensional situation. This leads to a direct proof of Theorem \ref{t1-1},
given in the next section, without using the three mathematical tools
just mentioned above.

\section{Non-linear case}

We now return to our general inequality (\ref{01}). First, we need
a measure $\hat\nu$, as in the last section, deduced from $\nu$.
Let ${\nu}^{\text{\#}}$ be the absolutely continuous part of $\nu$
with respect to the Lebesgue measure. Then, define
$${\hat\nu}(\d x)={\hat\nu}_p(\d x)= \bigg(\frac{\d {\nu}^{\text{\#}}}{\d x}\bigg)^{-1/(p-1)}\d x,\qqd p>1.$$
Next, we need a universal factor
$$k_{q, p}=\bigg[\frac{q-p}{p \text{\rm B}\!\Big(\frac{p}{q-p},\, \frac{p(q-1)}{q-p}\Big)}\Bigg]^{1/p- 1/q}
\le 2 \qqd \text{if $q\ge p$},$$
where ${\rm B}(\az, \bz)$ is the Beta function. In particular (as the limit of $q\downarrow p$),
$$k_{p, p}= p^{1/ p}\big(p^*\big)^{1/ {p^*}},$$
where $p^*$ is the conjugate of $p\in (1, \infty)$: $1/p+1/p^*=1$.

\thm\lb{t2-1}{\rm\rf{cmf13}{Theorem 2.6}}\qd {\cms Let $\mu (-M, N)<\infty$. Then the optimal constant $A$ in the Hardy-type inequality (\ref{01}) satisfies
\begin{itemize}
\item[(1)] the upper estimate $A \le k_{2, p} B^*$ for $1< p\le 2 \le q< \infty$ once
the pure point part of $\mu$ (denoted by $\mu_{\text{\rm pp}}$) vanishes, and
\item[(2)]  the lower estimate $A\ge B_*$ for $1<p, q<\infty$, where
$$\aligned
B^*&=
\sup_{x\le y}\frac{{\hat\nu}(x, y)^{(p-1)/p}}{\Big\{{\mu}(-M, x)^{\frac{p}{q(1-p)}}+{\mu} (y, N)^{\frac{p}{q(1-p)}}\Big\}^{(p-1)/p}},\\
B_*&=\sup_{x\le y}\frac{{\hat\nu}(x, y)^{(p-1)/p}}{\Big\{
{\mu}(-M, x)^{\frac{1}{1-q}}+{\mu} (y, N)^{\frac{1}{1-q}}\Big\}^{(q-1)/q}}.
\endaligned$$
\end{itemize}
Moreover, $B_*\le B^*\le 2^{1/p-1/q}B_*$ once $q\ge p$.
}
\dethm

Here are some remarks on Theorem \ref{t2-1}.
\begin{itemize}\setlength{\baselineskip}{0pt}
\item[(a)] The isoperimetric constants $B^*$ and $B_*$ are expressed
explicitly in measures $\mu$ and $\hat\nu$.
\item[(b)] The boundaries $-M$ and $N$ symmetric in the formulas of $B^*$ and $B_*$.
\item[(c)] Even through $B^*\ge B_*$ in general, but the rough ratio   $k_{q,p} 2^{1/p-1/q}$ of the upper and lower bounds is still $\le 2$.
\item[(d)] When $q=p$, we have
$$B^*=B_*=\sup_{x\le y}\frac{{\hat\nu}(x, y)^{(p-1)/p}}{\Big\{
{\mu}(-M, x)^{\frac{1}{1-p}}+{\mu} (y, N)^{\frac{1}{1-p}}\Big\}^{(p-1)/p}}.$$
 \item[(e)] Ignoring the $\mu(-M, x)$-term in the expression of $B^*$ or
  $B_*$, we obtain
  $$
  B^+=\sup_{y} {\hat\nu}(-M, y]^{1/ {p^*}}\,\mu (y, N)^{1/ q}.$$
Similarly, Ignoring the $\mu(y, N)$-term in the expression of $B^*$ or
  $B_*$, we obtain
  $$B^-=\sup_{x} {\hat\nu}(x, N)^{1/ {p^*}}\,\mu (-M, x)^{1/ q}.$$
  We have thus returned to one of the main results in the study of Hardy-type
  inequalities.
\end{itemize}

\thm\lb{t2-2}{\rm (1920---1992)}\qd {\cms For the Hardy-type inequalities
$$\|f\|_{L^q(\mu)}\le A^+ \big\|f'\big\|_{L^p(\nu)}, \qqd f(-M)=0$$
and
$$\|f\|_{L^q(\mu)}\le A^- \big\|f'\big\|_{L^p(\nu)}, \qqd f(N)=0,$$
we have the basic estimates $B^{\pm}\le A^{\pm}\le k_{q, p}B^{\pm}$.
Moreover, the factor $k_{q, p}$ is sharp.}
\dethm

There is a long history about the development of Theorem \ref{t2-2}. The reader is urged to refer to \ct{kmp} and \ct{cmf13} for a long list of
references including five books.

Having the experience in proving Theorem \ref{t1-1} and known the history of
Theorem \ref{t2-2}, it is hardly imaginable how to find a direct proof of
Theorem \ref{t1-1}, or even much more general Theorem \ref{t2-2}, without using capacity. To have a test, let us introduce the proof of a hard part
--- the upper estimate of Theorem \ref{t2-2}.

The idea is starting from Theorem \ref{t2-2}. For this, we split the interval
$(-M, N)$ into two parts: $(-M, \uz)$ and $(\uz, N)$, \vspace{-0.3truecm}
\begin{center}{\includegraphics[width=11.0cm,height=1.5cm]{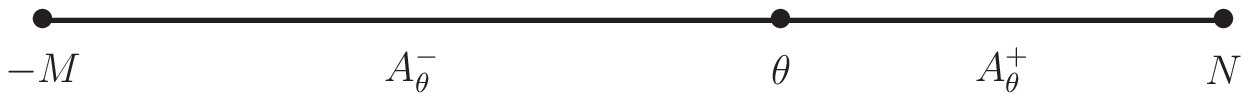}}\end{center}
Denote by $A_{\uz}^-$ the optimal constant on the left subinterval  $ (-M,\uz)$ and by  $A_{\uz}^+$ the one on the right subinterval $(\uz, N)$ with
 the same boundary condition $ f(\uz)=0$. Then, we can rewrite Theorem \ref{t2-2} as follows.
\medskip

\nnd {\bf Known Theorem}\qd {\cms Let $q\ge p$. Then  we have $B_{\uz}^{\pm}\le A_{\uz}^{\pm}\le k_{q, p} B_{\uz}^{\pm}$,  where
 $$\aligned
&  B_{\uz}^-=\sup_{x<\uz} {\hat\nu}(x, \uz)^{{1}/{p^*}}\mu(-M, x)^{1/q},\qqd
B_{\uz}^+=\sup_{y>\uz} {\hat\nu} (\uz, y)^{{1}/{p^*}}\mu(y, N)^{1/q}.\endaligned$$
}

\nnd{\bf Proof of the upper estimate}:
$k_{2,p}B^*\ge A$.

Rewrite $B^*$ as
$$\aligned
B^*&=
\sup_{x\le y}\frac{{\hat\nu}(x, y)^{(p-1)/p}}{\Big\{{\mu}(-M, x)^{\frac{p}{q(1-p)}}+{\mu} (y, N)^{\frac{p}{q(1-p)}}\Big\}^{(p-1)/p}}\\
&=\Bigg\{\sup_{x\le y}\frac{{\hat\nu}(x, y)}{{\mu}(-M, x)^{\frac{p}{q(1-p)}}+{\mu} (y, N)^{\frac{p}{q(1-p)}}}\Bigg\}^{(p-1)/p}\\
&=:\bigg\{\sup_{x\le y}\frac{{\hat\nu}(x, y)}{\fz(x)+\qz(y)}\bigg\}^{1/p^*}.
\endaligned$$
By proportional property, we have
$$\aligned
\frac{{\hat\nu}(x, y)}{\fz(x)+\qz(y)}
&=\frac{{\hat\nu}(x, \uz)+ {\hat\nu}(\uz, y)}{\fz(x)+\qz(y)}
\ge \bigg\{\frac{{\hat\nu}(x, \uz)}{\fz(x)}\bigwedge
     \frac{{\hat\nu}(\uz, y)}{\qz(y)}\bigg\},\qqd \uz\in (x, y).
\endaligned$$
Here $a\wedge b=\min\{a, b\}$
and similarly $a\vee b=\max\{a, b\}$.
Hence (omit what in $\{\cdots\}$)
$$\frac{{\hat\nu}(x, y)}{\fz(x)+\qz(y)}\ge \sup_{\uz\in (x, y)}\{\cdots\}.$$
Then
$$\aligned
\sup_{x\le y}\frac{{\hat\nu}(x, y)}{\fz(x)+\qz(y)}
&\ge \sup_{x\le y}\sup_{\uz\in (x, y)}\{\cdots\}\\
&=\sup_{\uz}\sup_{(x, y)\ni \uz}\{\cdots\}\\
&=\sup_{\uz}\bigg\{\bigg[\sup_{x\le\uz} \frac{{\hat\nu}(x, \uz)}{\fz(x)}\bigg]\bigwedge
\bigg[\sup_{y\ge\uz}  \frac{{\hat\nu}(\uz, y)}{\qz(y)}\bigg]\bigg\}.
\endaligned$$
Making power $1/p^*$ on both sides, by definition of $B_{\uz}^{\pm}$, we obtain
$$B^*\ge \sup_{\uz}\big(B_{\uz}^-\wedge B_{\uz}^+\big).$$
Here a problem appears: we need $\vee$ rather than $\wedge$ on the right-hand
side. To overcome this, we assume that
$\mu_{pp}=0$. Then, there exists $\bar\uz\in (-M, N)$ such that
$B_{\bar\uz}^-= B_{\bar\uz}^+$.
Therefore
$$B^*\ge \sup_{\uz}\big(B_{\uz}^-\wedge B_{\uz}^+\big)\ge B_{\bar\uz}^-
=B_{\bar\uz}^-\vee B_{\bar\uz}^+$$
and then
$$\aligned
k_{q, p}B^*
&\ge \big(k_{q, p} B_{\bar\uz}^-\big)\vee \big(k_{q, p} B_{\bar\uz}^+\big).\\
&\ge A_{\bar\uz}^- \vee  A_{\bar\uz}^+\qd\text{\rm(Known Theorem)}\\
&\ge \inf_{\uz} \big(A_{\uz}^-\vee A_{\uz}^+\big)\\
&\ge A\qd \text{\rm(splitting technique)}.\endaligned$$
Here each step holds for all $q\ge p$ except the last one. In which, some additional work is required, due to
the appearance of $f-\pi(f)$ rather than $f$ only. We prove the conclusion first for $q= 2\ge p$ and then extend it to $q\ge 2$, even to a large class of
normed linear space ${\mathbb B}$, as used in Theorem \ref{t1-3}, using a known lifting procedure (cf. \rf{cmf05}{\S 6.3}).  Note that in the proof above, we use a bridge $\uz$ to combine the known results on two subintervals together. But then remove it, otherwise, $\bar\uz$ for instance, may not be computable.
Nevertheless, it should be understandable that the present
analytic proof does not use the coupling, duality, or capacitary techniques.

Actually, much more topics are studied in \ct{cmf13}: the bilateral Dirichlet boun\-daries, logarithmic Sobolev inequalities, Nash inequalities, and so on.

\section{Improvements of the basic estimates}

Note that for two given numbers having smaller ratio, their difference can
be still quite big. Hence, it is meaningful to improve the basic estimates
introduced in the last two sections. In this section, we show the possibility
in doing so by a simplest example: $\mu=\d x$ and $\nu=\d x$ on $(0, 1)$.
We need to consider the following Hardy-type inequality
\be \|f\|_{L^q(\mu)}\le A \big\|f'\big\|_{L^p(\nu)}, \qqd f(0)=0 \lb{03}\de
only since the other cases (the ergodic case in particular) can be reduced to
this one by symmetry. The basic estimates for the optimal constant $A$ in (\ref{03}) are given in
Theorem \ref{t2-2}.

\xmp\qd{\cms Let $\mu=\d x$ and $\nu=\d x$ on $(0, 1)$. Then the optimal
constant $A$ in (\ref{03}) is given as follows.
\begin{itemize}
\item[(1)] When $p=q=2$, we have $A=2/\pi$.
\item[(2)] When $p=q\in (1, \infty)$, we have
   $$A=\frac{p}{\pi (p-1)^{1/p}}\sin\frac \pi p.$$
\item[(3)] For general $p, q\in (1, \infty)$, we have
$$A=\frac{p^{\frac{1}{q}} q^{1-\frac{1}{p}} (pq+p-q)^{\frac{1}{p}-\frac{1}{q}}
}
{(p-1)^{\frac 1 p}\, \text{\rm B}\left(\frac{1}{q},\,1-\frac{1}{p}\right)
 }.$$
\end{itemize}
}
\dexmp
This simplest example already shows that it is nontrivial from the special case
$p=q=2$ to the general one.

\prp\lb{t3-2}{\rm \ct{cmf13b}}\qd {\cms For the optimal constant $A$ in the last example, we
have the following improved estimates:
$$B\le {\bar\dz}_1\le A\le A^*\le \dz_1\le k_{q,p}B,$$
where
$$\aligned
B&=\frac{p^{1/q}((p-1)q)^{1-1/p}}{(pq+p-q)^{1-1/p+1/q}},\\
{\bar\dz}_1&=\frac{p^{1/q}((p-1)(q+1))^{1-1/p}}{(pq+p-q)^{1-1/p+1/q}},\\
A^*&=\bigg[\frac{p^*}{q}\bigg]^{\frac 1 q}
\bigg[\frac{p^*+q}{\pi p^*}\sin \frac{\pi p^*}{p^*+q}\bigg]^{\frac{1}{p^*}+ \frac 1 q}\qd\fbox{$=A$ if $q=p$}\\
\dz_1 &=\frac{1}{({q\gz^*}/{p^*}+1)^{\frac 1 q}}\bigg[\!\sup_{x\in (0, 1)}\frac{1}{x^{\gz^*}}\int_0^x\big(1-y^{{q/\gz^*}{p^*}+1}\big)^{\frac{p^*}{q}}\d y\bigg]^{\frac{1}{p^*}},
\qd \gz^*:=\frac{q}{p^*+q}
\endaligned$$
}
\deprp

The results in Proposition \ref{t3-2} are shown by Figures 1--4.

\begin{center}{\includegraphics[width=11.0cm,height=7.5cm]{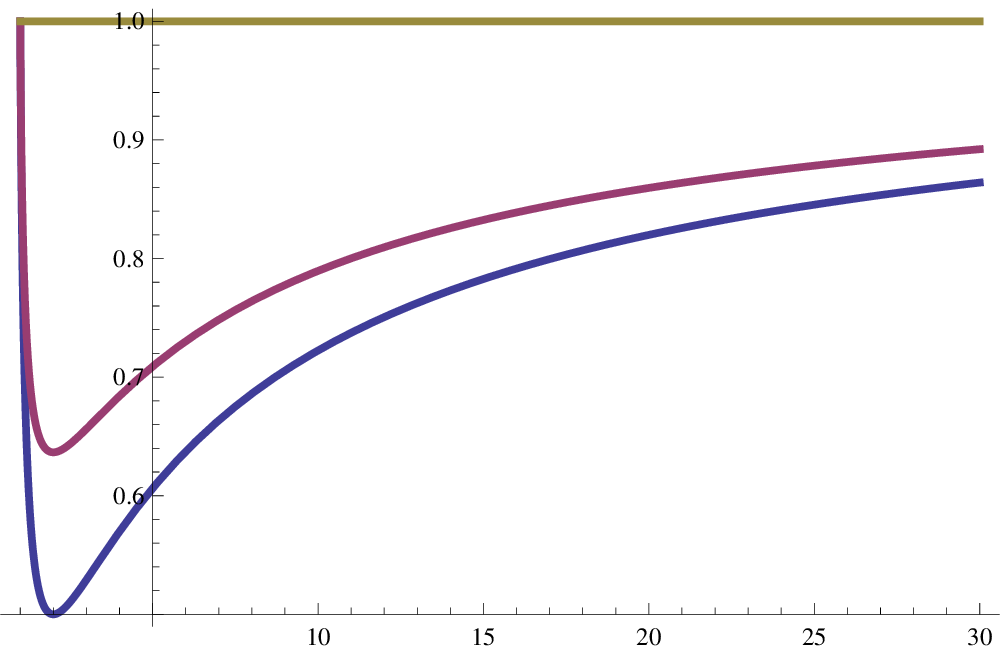}\newline
{\bf Figure 1}\qd The basic estimates of $A$: $p=q\in (1, 30)$.}\end{center}

Note that in the case that $p=q$, we have $A^*=A$. So there are five curves only in Figure 2.

The improvements are surprisingly effective.
Note that a suitable convex mean of the new upper and lower bounds should provides a quite precise appro-

\newpage

\nnd ximation of $A$. However, the convex means of the basic estimates do not have this property. When $p=q$, much more refined results can be found from \ct{cwz13a, cwz13b, cwz13c}.

\begin{center}{\includegraphics[width=11.0cm,height=7.5cm]{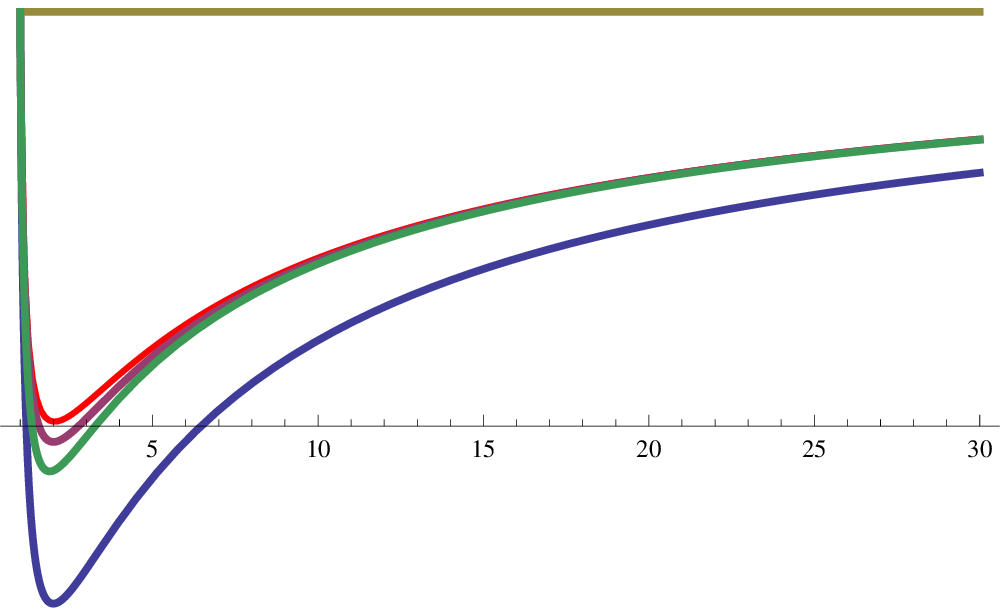}\newline
{\bf Figure 2}\qd The basic estimates of $A$ and their improvements: $p=q\in (1, 30)$. }\end{center}

Next, since $q\ge p$, we may write $q=p+r$ for some $r\ge 0$. In Figures 3 and 4, there are six curves, three of them are upper estimates and two of them are lower ones. The third curve from the bottom is the exact one; the top and the bottom curves consist of the basic estimates of the exact one. The other three curves are the improvements of the basic estimates.

\begin{center}{\includegraphics[width=11.0cm,height=7.5cm]{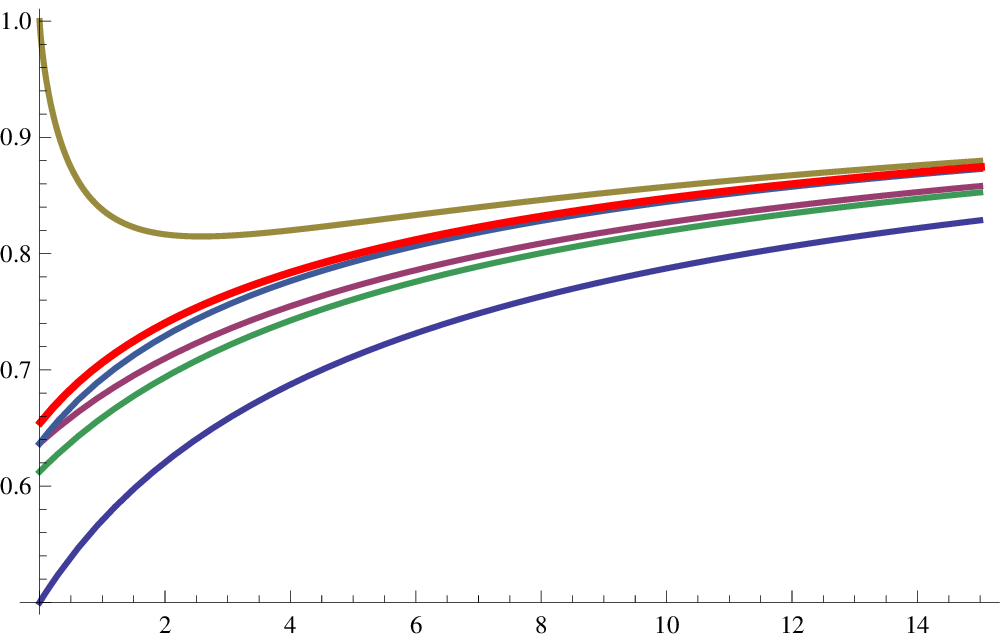}\newline
{\bf Figure 3}\qd The basic estimates and their improvements: $p=2$, $r\in (0, 15)$. }\end{center}

\begin{center}{\includegraphics[width=11.0cm,height=7.5cm]{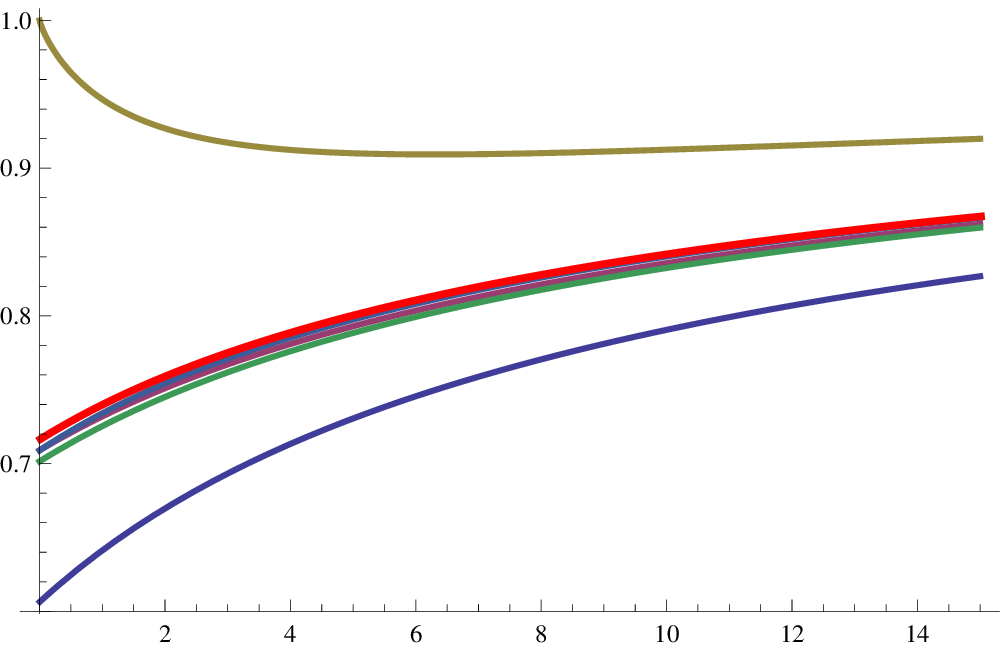}\newline
{\bf Figure 4}\qd The basic estimates and their improvements: $p=5$, $r\in (0, 15)$. }\end{center}

Note that in Figure 4, the new upper bounds and lower bounds are almost overlapped with $A$. In general, they are closer when $p$ or $q\ge p$ is larger.

Finally, we mention that the main results in this note: Theorems \ref{t1-1} and \ref{t2-1}, and Proposition \ref{t3-2} are new addition to the context of Hardy-type inequalities.
\bigskip

\nnd{\bf Acknowledgments}. {\small
This short note introduces an example to illustrate the great help we have obtained from Professor Masatoshi Fukushima in the past decades. It is a nice
chance to express our deepest thanks to him.
Research supported in part by the
         National Natural Science Foundation of China (No. 11131003),
         and by the ``985'' project from the Ministry of Education in China.
}

\nnd {\small School of Mathematical Sciences, Beijing Normal University,\newline
Laboratory of Mathematics and Complex Systems (Beijing Normal University),\newline
\mbox{\qqd} Ministry of Education, Beijing 100875,
    The People's Republic of China.\newline E-mail: mfchen@bnu.edu.cn\newline Home page:
    http://math.bnu.edu.cn/\~{}chenmf/main$\_$eng.htm
}


\begin{thebibliography}{000}
\setlength{\itemsep}{-0.8ex}
{\small

\bibitem{cmf04}
Chen, M.F. (2004).
{\it From Markov Chains To Non-Equilibrium Particle Systems}, 2$^{\rm nd}$ ed. World Sci., Singapore.\lb{cmf04}

\bibitem{cmf05}
Chen, M.F. (2005a).
{\it Eigenvalues, Inequalities, and Ergodic Theory}.
    Springer, London.\lb{cmf05}

\bibitem{cmf05b}
Chen, M.F. (2005b).
{\it Capacitary criteria for Poincar\'e-type inequalities}.
Potential Analysis 23(4), 303-322.\lb{cmf05b}

\bibitem{cmf10}
Chen, M.F. (2010).
{\it Speed of stability for birth--death processes.}
Front. Math. China 5(3), 379--515.\lb{cmf10}

\bibitem{cmf12}
Chen, M.F. (2012).
{\it Basic estimates of stability rate for one-dimensional diffusions.}
Chapter 6 in ``Probability Approximations and Beyond'', 75--99, Lecture Notes in Statistics 205,
eds. A.D. Barbour, H.P. Chan and D. Siegmund.\lb{cmf12}

\bibitem{cmf13}
Chen, M.F. (2013a).
{\it Bilateral Hardy-type inequalities}.
Acta Math. Sin. Eng. Ser. 2013, 29:1, 1--32.
\lb{cmf13}

\bibitem{cmf13b}
Chen, M.F. (2013b).
{\it The optimal constant in Hardy-type inequalities},
  in preparation.
\lb{cmf13b}

\bibitem{cmwf97}
Chen, M.F. and Wang, F.Y. (1997).
Estimation of spectral gap for elliptic operators.
Trans. Amer. Math. Soc. 349(3), 1239--1267.\lb{cmwf97}

\bibitem{cwz13a}
Chen  M.F., Wang L.D., Zhang Y.H. (2013a)
{\it Mixed principal eigenvalues in dimension one}.
Front. Math. China, 2013, 8(2): 317-343.\label{cwz13a}

\bibitem{cwz13b}
Chen, M.F., Wang, L.D., Zhang, Y.H. (2013b)
{\it Mixed  eigenvalues of discrete {$p$\,}-Laplacian},
preprint. \lb{cwz13b}

\bibitem{cwz13c}
Chen, M.F., Wang, L.D., Zhang, Y.H. (2013c)
{\it Mixed  Eigenvalues of  {$p$\,}-Laplacian},
preprint. \lb{cwz13c}

\bibitem{fm80}
Fukushima, M. (1980).
{\it Dirichlet Forms and Markov Processes}.
North Holland Math. 23.\lb{fm80}

\bibitem{fot11}
Fukushima, M., Oshima, Y. and Takeda, M. (2011).
{\it Dirichlet Forms and Symmetric Markov Processes,} 2$^{\text{nd}}$ ed.
Walter de Gruyter.\lb{fot11}

\bibitem{fmut}
Fukushima, M. and Uemura, T. (2003).
{\it Capacitary bounds of measures and ultracontractivi\-ty of time changed processes.}
J. Math. Pure et Appliqu\'ees 82(5), 553--572.\lb{fmut}

\bibitem{kmp}
Kufner, A., Maligranda, L. and Persson, L.E. (2007).
{\it The Hardy Inequality: About its History and Some Related Results}.
Vydavatelsky Servis.\lb{kmp}
}
\end{thebibliography}
\end{document}